\documentclass{article}
\usepackage{amsmath,amsfonts,amssymb}

\numberwithin{equation}{section}

\newcommand{\sezione}[1]{\section{\normalsize #1}}

\newtheorem{theorem}[equation]{\textnormal{THEOREM}}
\newtheorem{lemma}[equation]{\textnormal{LEMMA}}
\newtheorem{corollary}[equation]{\textnormal{COROLLARY}}
\newtheorem{proposition}[equation]{\textnormal{PROPOSITION}}
\newtheorem{definition}[equation]{\textnormal{DEFINITION}}
\newtheorem{example}[equation]{\textnormal{EXAMPLE}}
\newtheorem{remark}[equation]{\textnormal{REMARK}}
\newtheorem{remarks}[equation]{\textnormal{REMARKS}}
\newtheorem{note}[equation]{\textnormal{NOTE}}

\newcommand{\codim}{\operatorname{codim}}

\newcommand{\di}{\operatorname{div}}
\newcommand{\Pic}{\operatorname{Pic}}
\newcommand{\Di}{\operatorname{Div}}
\newcommand{\Cl}{\operatorname{Cl}}
\newcommand{\CaCl}{\operatorname{CaCl}}

\newcommand{\ZZ}{{\mathbb Z}}
\newcommand{\PP}{{\mathbb P}}
\newcommand{\QQ}{{\mathbb Q}}

\renewcommand{\o}{\omega}

\renewcommand{\O}{{\mathcal O}}
\newcommand{\I}{{\mathcal I}}
\newcommand{\F}{{\mathcal F}}
\newcommand{\E}{{\mathcal E}}
\newcommand{\G}{{\mathcal G}}
\newcommand{\Hom}{{\mathcal H}om}
\newcommand{\Coh}{{\mathcal C}oh}
\newcommand{\Refl}{{\mathcal R}efl}
\newcommand{\Div}{{\mathcal D}iv}

\newenvironment{pf}
{\noindent\textrm{Proof.}}
{\hfill{$\square$}\medskip}

\begin{document}

\vspace*{2.8cm}
\section*{Weil divisors on rational normal scrolls}
\begin{flushleft}
\medskip
RITA FERRARO 
Dipartimento di Matematica 
Universit\`a di Roma "Tor Vergata" 
Viale della Ricerca Scientifica, 00133 Roma. Con il supporto dell'Istituto
Nazionale Alta Matematica. 
\end{flushleft}
\bigskip
\bigskip
\sezione{{INTRODUCTION}}
Let
$X\subset \PP^n$ 
be a 
rational normal scroll
of degree $f$ and dimension $r=n-f+1$. $X$ is the image of
a projective bundle $\tilde X$ over $\PP^1$ of rank
$r-1$ through the birational  morphism $j$ defined by the tautological line
bundle
$\O_{\tilde X}(1)$. Depending on $\tilde X$, the scroll $X$ may be smooth
(in this case $j$ is an isomorphism) or singular.
The aim of
this note is to study Weil divisors on a {\it singular}
rational normal
scroll.
Let $V$ be the {\it vertex} of $X$,
  we have that $\codim(V, X)\geq 2$. When 
$\codim(V, X)=2$ the scroll $X$ is a {\it cone} over a rational normal
curve of
degree $f$; in this case let $E=j^{-1}(V)$ be the 
{\it exceptional divisor} 
in $\tilde X$.

The paper is divided into four sections. In section 2
we set up the notation and recall some standard facts 
about
Weil
divisors on normal varieties: the {\it strict image} map 
$j_\#:\CaCl(\tilde
X)\to \Cl(X)$ (defined in Prop \ref{j}),
 is used to describe the group $\Cl(X)$
 of  Weil divisors on $ X$ 
modulo linear equivalence.
 It turns out that there are two different cases:
(i) when $\codim(V, X)>2$, then $j_\#$ is an isomorphism;
(ii) when $\codim(V, X)=2$, then  $j_\#$ is surjective and
$\ker(j_\#)$ is
generated
by $E$. Since the  group $\CaCl(\tilde X)$ of Cartier divisors on $\tilde 
X$ modulo
linear equivalence is well known, this provides
 an explicit description of
$\Cl(X)$ (Cor. \ref{Cl}). Furthermore we introduce the basic theme of the
next section: the relation between the strict image, the {\it scheme
theoretic image} 
and the direct
image of the ideal sheaf 
of a divisor on $\tilde X$ (Lemma \ref{ideals}).

The goal of section 3 is to describe the sheaves corresponding to Weil
divisors on $X$.
It is known (see Prop. \ref{divsheaf}) that  
the group $\Cl(X)$ is in bijection with the
set $\Div(X)$
of
{\it divisorial
sheaves} on $X$ (i.e. coherent sheaves which are {\it reflexive} of
rank one), and induces on it a
natural group
structure. 
We   describe explicitly the group
$\Div(X)$ via the direct image morphism $j_*:\Pic(\tilde X)\to \Coh(X)$ of
sheaves. In particular, we will be able to say when the direct image 
$j_*\F$ of an invertible sheaf $\F\in\Pic(\tilde X)$  is {\it
reflexive}.
The analysis naturally splits in the two cases above mentioned.
When $\codim(V, X)>2$, 
we will prove with standard  techniques
 that $j_*\Pic(\tilde X)=\Div(X)$ (Cor. \ref{reflone}).
In particular (Prop. \ref{firstdiv}) 
the divisorial sheaf $\O_X(D)$
of a Weil divisor $D$ on $X$ is
\[ 
\O_X(D)\cong j_*\O_{\tilde X}(\tilde D),
\]
where $\tilde D$ is the {\it proper transform} of $D$ in $\tilde X$. 
When $\codim(V, X)=2$, this is no longer true. To overcome this problem we
 introduce the concept
of {\it integral total transform} $D^*\subset \tilde X$ of a Weil divisor
$D\subset X$ (see Def. \ref{deftot}).
We 
prove (Th. \ref{refltwo}) that $\Div(X)$ consists of the direct images
of those 
 line bundles $\F\in \Pic(\tilde X)$ such that the degree $\deg(\F_{|E})$
of $\F$
restricted to the
exceptional divisor $E$ is $<f$, where $f$ is the degree of the scroll
$X$.
In particular we will prove the projection formula:
\[
\O_X(D)\cong j_* (\O_{\tilde X}(D^*)). 
\]

In section 4 we will study the intersection of Weil divisors on $X$ in
the critical case, i.e. 
$\codim(V, X)=2$.
According to 
\cite{h2}, Prop. 2.4.
an effective Weil divisor $D$ on $X\subset\PP^n$ 
is itself a closed
subscheme of $X$ of pure codimension $1$ with no embedded components.
 Given two effective divisors
$D$ and $D^\prime$
on $X$ with no common components, we can consider the 
scheme-theoretic intersection $Y=D\cap D^\prime\subset\PP^n$ and ask, for
example, for
its {\it degree}. Note that
 when $\codim(V, X)>2$ one can compute the degree of a "complete
intersection" of $l$ ($1\leq l\leq r-1$) divisors $D_1, \dots , D_l $
using the natural intersection form on $\Cl(X)$ inherited from $\tilde
0X$, via the isomorphism $j_\#$. When
$\codim(V, X)=2$, the unique linear theory of intersection which 
can be defined in $X$ is the generalization to higher dimension
of the theory developed by Mumford in \cite{m} for a normal surface.
According to this theory, the intersection number of
two divisors
$D, D^\prime$ on $X$ is a rational number defined as
the intersection number 
 of the corresponding
{\it rational total transforms}  in $\tilde X$ (see Def. \ref{mtot}),
which
is in general a rational number and does not represent
the degree of the  scheme theoretic intersection $Y=D\cap D^\prime$. 
In Th. \ref{reso} we will use the integral total transform
to find the minimal
{\it reflexive resolution} of $\O_Y$
as an $\O_X$-module,
which allows us  to compute the degree of $Y$ (Prop. \ref{grado}).

The last section is devoted to examples and applications.
 In particular
in Ex.\ref{primo} we show that every  divisor of degree $\geq n$ on a
rational normal cone
of degree $n-1$ in $\PP^n$ has maximal arithmetic genus.
In Ex. \ref{secondo} we show that, when $\codim(V, X)=2$, every
effective divisor
$D\subset X$ and every scheme theoretic intersection $Y=D\cap
D^\prime\subset X$
of two effective divisors with no common components
are arithmetically Cohen-Macaulay schemes in $\PP^n$.
In Ex. \ref{quarto} we will compute
the arithmetic genus  of $Y$.

Much of this material was motivated 
by the subject of my doctoral
thesis: 
{\it Classification of curves of maximal genus in
$\PP^5$}, since these curves lie on (possibly singular) rational normal
three-folds. In that context it is necessary to compute the degree of the
scheme theoretic intersection of two divisors on a rational normal cone
$X$
over a twisted cubic in $\PP^5$,  and to develope some
linkage tecniques 
on $X$. Moreover, for the linkage problem,  it is necessary to know the
divisorial sheaves on $X$.
Thanks to my advisor
Ciro Ciliberto.

\sezione{PRELIMINARIES} 

In this section we recall some basic facts about rational 
normal scrolls and we describe the group $\operatorname{Cl}(X)$
in terms of $\operatorname{CaCl}(\tilde X)$.
For more details about rational normal scrolls the reader may consult
for example \cite{r} or \cite{eh}.  

A rational normal scroll $X\subset \PP^n$ 
is the  image of a projective bundle
$ \pi:\PP({\cal E})\to \PP^1$ over $\PP^1$ of rank $r-1$
through the morphism $j$ defined by
the tautological bundle
${\O}_{\PP ({\cal E})}(1)$, where
$\E\cong \O_{\PP^1}(a_1)\oplus \cdots \oplus
\O_{\PP^1}(a_{r})$
with $0\leq a_1\leq \cdots \leq a_{r}$ and $\sum 
a_i=f$. 
If $a_1=\cdots=a_l=0$, $1\leq l <r$, $X$ is singular and
the  vertex
$V$ of $X$ has dimension $l-1$. Let us denote $\PP(\E)=\tilde X$.
The morphism  $j: \tilde{X}\to X$ is a rational
resolution of
singularities, i.e. $X$ is {\it normal} and 
{\it arithmetically Cohen-Macaulay}
and $R^ij_*\O_{\tilde X}=0$ for $j>0$. We will call $j:\tilde X\to X$
the {\it canonical resolution} of $X$.

It is a general fact that Weil divisors on a normal scheme
do not depend on closed subsets of codimension $\geq 2$. We refer to
\cite{h1}, 
Prop. II, 6.5. for this basic fact which we will  continuously
use in this note. In Th. \ref{refl} we will see the equivalent result
 in terms
of sheaves. Let $X_S=X\setminus V$ be the smooth part of $X$;
since $X$ is normal by \cite{h1}, Prop. II, 6.5.
there is an isomorphism $Cl(X)\to Cl(X_S)$ defined by
$D=\sum n_iD_i\to\sum n_i (D_i\cap X_S )$ , where
$D_i$ are prime divisors on $X$; 
in other words, given a prime divisor $D_S$ 
on $X_S$ there is an unique prime divisor
$D=\overline{D_S}$ (the closure of $D_S$ in $X$ with the induced scheme
structure), which extends $D_S$ on all $X$.  
Since $j_{|j^{-1}X_S}$ is an isomorphism, we can rewrite
\cite{h1}, Prop. II, 6.5.
in the following way:

\begin{proposition}
\label{j}
Let $X$ be a singular rational normal scroll, let $V$ be its vertex and 
let $X_S$ be its smooth part. 
Let $j:\tilde X\to X$ be the canonical resolution. 
Then:
\begin{enumerate}
\item there is a surjective homomorphism
$j_\#:\CaCl(\tilde X)  \to 
\operatorname{Cl}(X)$
defined by $C=\sum n_i C_i\to \sum n_i 
\overline{j ({C_i}\cap{j^{-1}X_S})}$, where
we ignore those $C_i\cap j^{-1}X_S$ which are empty;
\item if $\codim(V, X)>2$, then $j_\#:\Pic(\tilde X)\to \Cl(X)$
is an isomorphism;
\item if $\codim(V, X)=2$ and $E$ is the exceptional divisor
of $j$, then there is an exact sequence:
$
0\to\ZZ\to\CaCl(\tilde X)\stackrel{{j_\#}}\to \Cl(X)\to 0
$ 
where the first map is defined by $1\mapsto 1\cdot E$.
\end{enumerate}
\end{proposition}

Given a Cartier divisor $C$ on $\tilde X$,  the Weil divisor
$j_\#(C)$ on $X$ will be called the {\it strict image} of $C$ through $j$. 
It is well known that
$\Pic(\tilde{X})=\mathbb{Z}[\tilde{H}]\oplus\mathbb{Z}[\tilde{R}]$,
where $[\tilde{H}]=[\O_{\tilde{X}}(1)]$ is the hyperplane class and
$[\tilde{R}]=[\pi^*\O_{\PP^1}(1)]$ is the class of the fibre of
the  map $\pi: \tilde X\to \PP^1$.
Let $H=j_\#\tilde H$ and $R=j\#\tilde R$ 
be the strict images of $\tilde H$  and $\tilde R$ respectively
(i.e. respectively an hyperplane section 
and a divisor in the 
ruling of $X$). Then as a consequence of Prop. \ref{j} we have the
following:

\begin{corollary}
\label{Cl}
Let $X\subset \PP^n$ be a singular rational normal scroll of
degree $f$  and let $j:\tilde X\to X$
be its canonical resolution. Then
\begin{enumerate}
\item
if $\codim (V, X)>2$,
$\operatorname{Cl}(X)\cong\ZZ [H]\oplus\ZZ [R]$;
\item
if $\codim(V, X)=2$, $E\sim\tilde H-f\tilde R$ and
$\operatorname{Cl}(X)\cong\mathbb{Z}[R]$.
\end{enumerate}
\end{corollary}
\begin{pf}
1) Follows immediately from Prop. \ref{j} 2).
 When $\codim(V, X)=2$ an hyperplane section $H$ passing
trough $V$ splits in
the union
of $f$ fibers $R$, therefore we have $H\sim fR$, i.e. by Prop. \ref{j}
3) $E\sim \tilde H-f\tilde R$.
\end{pf}
\begin{remark}
\label{nocartier}
The Weil divisor $R$ on $X$ is not Cartier, since it is not locally
principal
in a neighboorhood of the vertex $V$. 
The proof is  similar to the classical example of a quadric cone
(\cite{h1}, Ex. II, 6.5.2).
\end{remark}

However, besides the strict image, there are other natural ways of
"pushing down" a divisor on $\tilde X$. Namely one can consider the
scheme theoretic image.
If $C\subset \tilde X$ is a   subscheme of 
$\tilde X$, 
the  {\it scheme-theoretic image}  $j_*(C)$
 of $C$ through $j$  is the unique
closed subscheme 
of $X$ with the following property: the morphism $j\circ i$, where
$i: C\hookrightarrow \tilde X$ is the canonical injection of $C$ in
$\tilde X$,
factors trough $j_*(C)$, and if $D$ is any other closed subscheme
of $X$ through which $j\circ i$ factors, then $j_*(C)\hookrightarrow X$
factors through $D$ also (see \cite{ega}, \S 6.10 pp. 324-325 and 
\cite{h1} Ex. II, 3.11 d)). 
According to \cite{h2}, Prop. 2.4.
an effective Weil divisor $D$ on $X\subset\PP^n$, i.e. 
a formal sum $\sum n_i D_i$
with $n_i\in \ZZ^+$ and $D_i$ irreducible, corresponds to a closed
subscheme of $X$ of pure codimension $1$ with no embedded components.
Therefore we can talk about the scheme theoretic image $j_*(C)$ of an
effective
divisor $C$.

Given a prime  divisor $D$ on $X$, the {\it proper transform}
$\tilde D$
of $D$ on $\tilde X$ is the 
closure $\overline{j^{-1}(D\cap X_S)}$
with the induced reduced scheme structure. The proper transform
of any Weil divisor in $X$ is then defined by linearity.
If $\tilde D$ is effective it follows that $\tilde D$ is the
{\it scheme-theoretic 
closure} of ${\tilde D}_{|j^{-1}X_S}$ in $\tilde X$, that is the smallest
closed subscheme of $\tilde X$ containing ${\tilde D}_{|j^{-1}X_S}$
as subscheme  (see \cite{ega}
Def. 6.10.2 pg. 324); this implies that $\tilde D$ is the scheme-theoretic
image
of ${\tilde D}_{|j^{-1}X_S}$ through the canonical injection $k:
{\tilde D}_{|j^{-1}X_S}\hookrightarrow \tilde X$.

\begin{remarks}
\label{image}
\begin{enumerate}
\item
For any Weil divisor $D$ on $X$, 
 from the definition
of $j_\#$ and of $\tilde D$ it follows that
 $j_\#(\tilde D)=D$.
\item
Let $\codim(V, X)>2$ and let $D$ be a divisor on $X$.
 The unique Cartier divisor $C$ on $\tilde X$
such that $j_\#(C)=D$ 
is the proper transform $\tilde D$.
Indeed, since $C$ and $\tilde D$ coincide outside a
subset of codimension
$\geq 2$, they must be equal.
\item If $\codim(V, X)=2$ we have 
$j_\#(\tilde D +nE)=D$
for every $n\in\ZZ$. In fact $E\cap j^{-1}X_S=\emptyset$.
\item
If $D$ is  effective, then $j_*(\tilde D)=D$. 
Let us denote by $k: {\tilde D}\cap{j^{-1}X_S}\hookrightarrow
\tilde X$
and $i: {D}\cap{X_S}\hookrightarrow X$ the canonical injections and by 
$h:\tilde D\cap j^{-1}X_S\to D\cap X_S$ the isomorphism
$j_{|\tilde D\cap j^{-1}X_S}$;  
we then have $j\circ k=i\circ h$. By \cite{ega} Prop. 6.10.3
pg.324 (transitivity of the scheme-theoretic images) we  have that
$j_*(k_*({\tilde D}\cap{j^{-1}X_S}))=
i_*(h_*({\tilde D}\cap{j^{-1}X_S}))$. Since,  as we noted before,
$k_*( {\tilde D}\cap{j^{-1}X_S})=\tilde D$ and
(for the same reason) $i_*({D}\cap{X_S})=D$, this proves what we claimed. 
\end{enumerate}
\end{remarks}
In the following Lemma we prove
that for an arbitrary divisor $C$ on $\tilde X$ the
 scheme theoretic image $j_*(C)$ 
is defined by the
ideal sheaf $j_*\I_{C|{\tilde X}}\hookrightarrow j_*\O_{\tilde X}=\O_X$.
As it will be more clear in the next section,   
the strict image $j_\#(C)$ is defined by
${j_*\I_{C|{\tilde X}}}^{\vee\vee}$ (Prop. \ref{divsheaf} 4)).

\begin{lemma}
\label{ideals}
Let $C$ be an effective Cartier divisor on $\tilde X$ and let $j_*(C)$ be
its
scheme-theoretic image in $X$.
Then $j_*(\I_{C|{\tilde
X}})\cong\I_{j_*(C)|X}$.
\end{lemma}
\begin{pf}
Let us consider the following diagram of sheaves on $X$:
\[
\begin{matrix}
0\to & \I_{j_*(C)|X} & \to &\O_X &\to & \O_{j_*(C)} & \to 0 \\
     & \downarrow{j^\#} & & \downarrow{j^\#} && \downarrow{{j^\#}_{|D}} & \\
0\to & j_*\I_{C|{\tilde X}} & \to & j_*\O_{\tilde X} &\to & j_*\O_{C} &
\to
\end{matrix}
\]
Since $j_*\O_{\tilde X}\cong \O_X$, the morphism $j^\#: \O_X \to
j_*\O_{\tilde X}$ 
is
an isomorphism; therefore,
by the Snake's Lemma, we
have to prove that
the morphism $\O_{j_*(C)}\to j_*\O_C$ is injective.
This follows because otherwise
$\ker(j^\#_{|j_*(C)})\hookrightarrow \O_{j_*(C)}$ would define a
subscheme
$C^\prime$ of $j_*(C)$ such that 
the morphism $\O_{j_*(C)}\to j_*\O_C$ factors in
$\O_{j_*(C)}\to \O_{C^\prime} \to j_*\O_C$, but this cannot happen by
universal property of $j_*(C)$.
\end{pf}

\sezione
{DIVISORIAL SHEAVES}

We consider here the problem of describing the group $\Div(X)$
of divisorial sheaves
on a singular rational normal scroll $X$ in terms
of the Picard group $\Pic(\tilde X)$ of the canonical resolution $\tilde 
X$. It is known (see
Prop. \ref{divsheaf} below) that $\Div(X)$ is naturally isomorphic
to the group $\Cl(X)$ 
of Weil divisors modulo linear equivalence,
The analysis naturally splits in two cases: $\codim(V, X)>2$ and
$\codim(V, X)=2$. First we deal with the first case, which can be treated
in a more or less standard way. The result (Cor. \ref{reflone}) is that
the
natural map $j_*:\Pic(\tilde X)\to \Coh(X)$ is in fact an isomorphism   
$j_*:\Pic(\tilde X)\to \Div(X)$. This is no longer true when
$\codim(V, X)=2$. To overcome this problem we will introduce the concept
of integral total transform of a Weil divisor $D$ on $X$.
But first let use give  basic definitions and properties of
divisorial sheaves on a normal scheme; for details and for a more
general point of view
the reader may consult \cite{h2}, \S 2.
Let $X$ be a normal scheme. We recall that
a coherent sheaf $\F$ on  $X$ is {\it reflexive} 
if the natural map $\F\to\F^{\vee\vee}$ is an isomorphism, where ${\F^\vee}$ 
denotes the dual sheaf $\Hom(\F,
\O_X)$.
The following Theorem, which says that reflexive sheaves depend
only on subsets of codimension $1$, is a basic fact which we will
use in this section. 
\begin{theorem}
\label{refl}
Let $X$ be a normal scheme and let $Y\subset X$ be a closed subset
of codimension $\geq 2$. Then the restriction map induces an
equivalence of categories from the category $\Refl(X)$
of reflexive sheaves on $X$ to the category 
$\Refl(X\setminus Y)$
of reflexive sheaves on $X\setminus Y$.  
\end{theorem}
\begin{pf}
\cite{h2}, Th. 1.12.
We recall from this proof the way to extend a reflexive 
sheaf $\G$ on $X\setminus Y$
to a reflexive sheaf $\F$ on $X$. It consists in taking a coherent extension
$\F_0$ of $\G$ in $X$ (which exists by a general result on extensions 
of coherent sheaves), and then to put $\F={\F_0}^{\vee\vee}$.
The double dual of any coherent sheaf
is in fact always reflexive.
\end{pf}
\begin{definition}
Let $X$ be a normal scheme.
Let $D$ be a Weil divisor on $X$.
If $K(X)$ denotes the function field of $X$ (\cite{h1}, pg. 91 and pg. 141),
then the sheaf $\O_X(D)$ defined for  every open set $U\subset X$ as 
\[
\Gamma(U, \O_X(D))=\{{f\in K(X)| \di f +D\geq 0 \quad
 \text{on}\quad  U}\}.
\]
is called the divisorial sheaf of $X$.
\end{definition}
The following Proposition
(see \cite{h2}, Prop. 2.8.)
 describes the equivalence between 
reflexive sheaves and divisorial sheaves.
In point 4) is defined the {\it group structure}
on $\Div(X)$ induced by $\Cl(X)$.

\begin{proposition}
\label{divsheaf}
Let $X$ be a normal scheme.
\begin{enumerate}
\item
For any Weil divisor $D$ the sheaf $\O_X(D)$
 is reflexive and locally free of rank one
at every generic point and at every point of codimension 1.
\item Conversely, every reflexive sheaf which is locally free
of rank one
at every generic point and at every point of codimension 1 is isomorphic
to  $\O_X(D)$ for some Weil divisor $D$.
\item
If $D_1$ and $D_2$ are Weil divisor on $X$, $D_1\sim D_2$
if and only if $\O_X(D_1)\cong\O_X(D_2)$ as $\O_X$-modules.
\item
If $D$, $D_1$, $D_2$ are Weil divisors on $X$,
 $\O_X(-D)\cong{\O_X(D)}^\vee$ and
$\O_X(D_1+D_2)={(\O_X(D_1)\otimes\O_X(D_2))}^{\vee\vee}$.
\end{enumerate}
\end{proposition}

We come back now to the case of a singular rational normal scroll
$X$. Our goal is to explicitly describe $\Div(X)$ with its group
structure and the idea is to compare $\Div(X)$ with $\Pic(\tilde X)$, via
the direct image map
$j_*$ of sheaves.
 There is a natural surjective map, induced by the
strict image map
 $j_\#:\Pic(\tilde X)\to Cl(X)$,
which we call
again $j_\#$:
\begin{eqnarray*}
j_\#:\Pic(\tilde X)&\to &\Div(X)\\
\O_{\tilde X}(C)& \mapsto &{j_*(\O_{\tilde X}(C))}^{\vee\vee}
\end{eqnarray*}
which is  bijective when $\codim(V, X)>2$.
The sheaf ${j_*(\O_{\tilde X}(C))}^{\vee\vee}$ must be 
the divisorial sheaf $\O_X(j_\#(C))$ by Th. \ref{refl} 
since it 
is reflexive and  isomorphic to $\O_X(j_\#(C))$ outside a subset
 of codimension $\geq 2$
(in the open set $X_S$). We are then going to compare
$j_\#:\Pic(\tilde X)\to \Div(X)$ with $j_*:\Pic(\tilde X)\to \Coh(X)$
or, in other words, we are going to check when $j_*(\O_{\tilde X}(C))$
is a reflexive sheaf. Of course we will have, according to Prop. \ref{j},
two separate cases: a) $\codim(V, X)>2$ and  b) $\codim(V, X)=2$.

First we briefly describe $j_*(\Pic(\tilde X))$. We refer for details
to \cite{s} or \cite{f}. 
Let us denote: 
\[
\O_{\tilde X}(a,b):=\O_{\tilde X}(a\tilde H+b\tilde R)
\]
 the invertible sheaf
associated to the class
$[a\tilde H+b\tilde R]$ in $\Pic(\tilde X)$ and
let us consider on $X$ their direct images: 
\[
\O_X(a, b):= j_*\O_{\tilde{X}}(a, b),
\]
with $a, b\in\mathbb{Z}$. 
The cohomology of 
$\O_{\tilde{X}}(a, b)$ can be explicitly  calculated using the Leray
spectral sequence, which,  since $R^i\pi_*\O_{\tilde X}(a, b)=0$
for every
$a, b\in \ZZ$ and $0<i<r-1$ (by Grauert's Theorem), simplifies as follows:
\begin{eqnarray*}
\dots & \to &  H^i(\pi_*\O_{\tilde X}(a, b))\to H^i (\O_{\tilde X}(a,
b))\to
H^{i-r+1}(R^{r-1}\pi_*\O_{\tilde X}(a, b))\to \\
&\to & H^{i+1}
(\pi_*\O_{\tilde X}(a, b)) \to \dots
\end{eqnarray*}
Therefore for $i<r-1$
and $a\geq 0$ we obtain:
\begin{equation}
\label{cohom}
H^i(\tilde X, \O_{\tilde X}(a, b))\cong
H^i(\PP^1, \pi_*\O_{\tilde X}(a, b))\cong
\sum_{|I|=a} H^i(\PP^1, \O_{\PP^1}( b +\sum_{j\in I}a_j))
\end{equation}
which are zero of course if $a< 0$ and for $1<i<r-1$; while for $j=r-1,
r$ we obtain:
\[
H^j(\tilde X, \O_{\tilde X}(a, b))\cong
H^{j-r+1}(\PP^1, R^{r-1}\pi_*\O_{\tilde X}(a, b)),
\]
which can be computed also by Serre duality using (\ref{cohom}).
For $a\geq 0$ and $b\geq -1$ using (\ref{cohom}) we 
compute:
\begin{equation}
\label{dim}
h^0(\tilde X, \O_{\tilde X}(a, b))=
f{{a+r-1}\choose{r}}+(b+1){{a+r-1}\choose{r-1}}.
\end{equation}
For $b<-1$, the dimension 
$h^0(\tilde X, \O_{\tilde X}(a, b))$  depends on the type of the
scroll,
i.e. on the integers $a_1, \dots, a_r$.
We recall from \cite{s} 
that we have the vanishing
\begin{equation}
\label{Ri}
R^i j_* \O_{\tilde X}(a, b)=0 
\end{equation}
for $i>0$ and for all  $a\in \mathbb{Z}$ and $b\geq -1$, 
which implies, via the degenerate Leray spectral sequence associated
to $j$:
\begin{equation}
\label{hi}
h^i(\O_X(a, b))=h^i(\O_{\tilde{X}}(a, b))
\end{equation}
for $i\geq 0$. Moreover in \cite{s} it is proved that the dualizing
sheaf $\o_X$ of $X$ is:
\begin{equation}
\label{dual}
\o_X=j_*\O_{\tilde X}(K_{\tilde X})=\O_X(-r, f-2).
\end{equation}
As we will see (Cor. \ref{reflone} and Th. \ref{refltwo}), $\o_X$ is
always
a divisorial
sheaf,  so it makes sense to talk about the {\it canonical divisor}
$K_X\sim -rH+(f-2)R$ on $X$.

Let us consider  the case when $\codim(V, X)>2$.
\begin{proposition}
\label{firstdiv}
Let $\codim(V, X)>2$, let $D\sim aH+bR$ be a Weil divisor on $X$ and let
$\tilde D\sim a\tilde H+b\tilde R$ be
its proper transform on $\tilde X$. Then 
\[
j_*\O_{\tilde X}(\tilde D)=\O_X(D)=\O_X(a, b).
\]
\end{proposition}
\begin{pf}
It is sufficient to consider the local situation. Let $U$ be an open set
containing the vertex $V$ and let $U^\prime=j^{-1}(U)$. Call
$V^\prime=j^{-1}(V)$.
Since $\codim(V^\prime, \tilde X)\geq 2$, then 
$H^0(U^\prime, \O_{\tilde X}(\tilde D))\cong H^0(U^\prime\setminus
V^\prime, \O_{\tilde X}(\tilde D))$; moreover
$H^0(U^\prime\setminus
V^\prime, \O_{\tilde X}(\tilde D))\cong 
H^0(U\setminus V, \O_{ X}( D))\cong H^0(U, \O_{ X}( D))$
where the last isomorphism follows from Th. \ref{refl}.
\end{pf}

\begin{corollary}
\label{reflone}
Let $\codim(V, X)>2$,
then $\O_X(a, b)$ is a reflexive sheaf for every
$a, b\in \ZZ$ and
\[
\Div(X)=\{\O_X(a, b)\, | \quad  a,b\in\ZZ\}.
\]
Moreover the natural group structure on $\Div (X)$ 
inherited from $\Cl(X)$
is given by
$
<\O_X(a, b), \O_X(a^\prime, b^\prime)>\mapsto
\O_X((a+a^\prime),(b+b^\prime))
$
and
$
{\O_X(a, b)}^\vee \cong \O_X(-a, -b).
$
\end{corollary}
\begin{pf}
The first assertions follows directly from Prop. \ref{firstdiv}.
The description of the group structure on $\Div(X)$ 
follows from Prop. \ref{firstdiv} and part 4) of Prop. \ref{divsheaf}.
\end{pf}

\begin{remark}
Note that $\O_X((a+a^\prime),(b+b^\prime))\cong{(\O_X(a,
b)\otimes\O_X(a^\prime,
b^\prime ))}^{\vee\vee}$, by
part 4) of Prop. \ref{divsheaf}. 
In general
$
\O_X(a, b)\otimes\O_X(a^\prime,b^\prime)
$ is not reflexive. For example $\O_X(0, 1)\otimes \O_X(0,
-1)$ is not reflexive. Otherwise we would have
$\O_X(0,
1)\otimes\O_X(0, -1)\cong {(\O_X(0,
1)\otimes\O_X(0, -1))}^{\vee\vee}\cong \O_X$. i.e. $\O_X(0, 1)$ would be 
invertible, which is a
contradiction since $R$ is not Cartier (see Remark \ref{nocartier}).
\end{remark}

When the codimension of $V$ is $2$ we know (Remarks \ref{image}  3))
that 
the map $j_\#:\Pic(\tilde X)
\to \Cl(X)$ is not injective and therefore also the map
$j_\#:\Pic(\tilde X)
\to \Div(X)$ is not injective. It turns out that the study of $\Div(X)$
becomes very simple if we introduce, for any Weil divisor
$D$ on $X$, a Cartier divisor $D^*$ on $\tilde X$, which we will call
{\it the integral total transform} of $D$, which plays, in some sense,
the role
of the proper transform in case $\codim(V, X)\geq 2$.

The problem of defining, for a Weil divisor $D$
trough $V$ on $X$, the total transform  $\tilde X$
has been considered and
solved by Mumford in \cite{m} on a normal surface, with the goal of
developing the bilinear intersection theory on normal surfaces
(see also \cite{sa}).
Mumford's theory can be generalized on rational normal cones
by defining the total transform $j^* D$ of a  divisor $D$ on $X$ as:
\begin{definition}
\label{mtot}
Let $\codim(V, X)=2$ and let $D$ be a Weil divisor on $X$. Then
the (rational) total transform of $D$ in $\tilde X$ is:
\[
j^* D=\tilde D + qE,
\]
where $E$ is the exceptional divisor on $\tilde X$ and $q$ is a
rational
number uniquely determined by the equation: $(\tilde D +qE)\cdot E\cdot 
{\tilde H}^{r-2}=0$.
\end{definition}
If $\tilde D\sim a\tilde H+b\tilde R$, then  we find that
$
q=\frac{b}{f}.
$
If $D$ is effective, since $\tilde D$ does not contain $E$, we have that
$b=\tilde D\cdot E\cdot {\tilde H}^{r-2}\geq 0$, i.e. $q\geq 0$.
Mumford's total transform $j^* D$
is in general  a $\QQ$-divisor and it is integral if and only if $D$
is Cartier. As we will see,  it is more convenient for our purposes
to consider,
if $D$ is effective, the round-up $\lceil j^*D\rceil$ of $j^*D$, i.e. the
smallest
integral divisor on $\tilde X$ containig $j^*D$. 
So let us give the following
definition:
\begin{definition}
\label{deftot}
Let $\codim (V, X)=2$.
Let $D$ be an effective Weil divisor on $X$,
we define the integral total transform $D^*$ of $D$
as:
\[
D^*=\lceil j^*D\rceil=
 \tilde D+ \lceil q\rceil E 
\]
where $\tilde D\sim a\tilde H+b\tilde R$ is the proper transform of $D$, 
$q=\frac{b}{f}$ is the same number appearing in Def.
\ref{mtot} and $\lceil q\rceil$ is the round-up of $q$, i.e. the smallest 
integer $\geq q$.
We define the total transform  of $-D$ as
${(-D)}^* =-D^*$.
\end{definition}

\begin{note}
It is a simple computation to get the following
equivalent expression of $D^*$. Let $D\sim dR$
be effective, i.e. $d\geq 0$, and divide
$d-1=kf+h$ ($k\geq -1$ and $0\leq h<f$).
Then:
\begin{equation}
\label{secdef}
D^*\sim (k+1)\tilde H-(f-h-1)\tilde R.
\end{equation}
The relations between these coefficients  and the ones in  
Def. \ref{deftot} are  $k+1=a+\lceil q\rceil$ and
$f-h-1=f\lceil q\rceil -b$.
From formula (\ref{secdef}) it is clear that $D^*$ is uniquely determined
by 
the  class of linear equivalence of $D\sim dR$, i.e. by the degree $d$.
\end{note}

\begin{proposition}
\label{totale}
Let $\codim (V, X)=2$. Let $D$ be an effective Weil divisor on $X$. Then
the integral total transform $D^*=\tilde D+\lceil q\rceil E$ 
is the biggest Cartier divisor $C$ on $X$ such that $j_*(C)=j_\#(C)=D$.
More precisely:
\[
j_*(\tilde D+\alpha E)=j_\#(\tilde D+\alpha E)=D
\] 
if and only if $0\leq \alpha \leq \lceil q\rceil$; in this
case 
$
j_*\I_{\tilde D+\alpha E |\tilde X}=\I_{D|X}.
$
\end{proposition}
\begin{pf}
In Remarks \ref{image}  we have seen that
$j_\#(\tilde D)=j_*(\tilde D)=D$ and 
$j_\# (\tilde D +mE)=D$ for every $m\in \ZZ$. Therefore
it is enough to prove that
$j_*(D^*)=D$ and that $D$  is a proper subscheme of $j_*(D^* + mE)$  
for every integer $m>0$.
Proving the first equality is equivalent to prove that
$j_*(D^*)$ does not have embedded components. 
Let us fix a divisor $D^\prime\sim (f-h-1)R$
on $X$; then $F=D+D^\prime\sim (k+1)H$ is a Cartier divisor on $X$
cut out by a hypersurface $F$ of degree $k+1$. 
Let us consider on $\tilde X$
the divisor $j^*F\sim (k+1)\tilde H$.
We have that $D^* +{\tilde D}^\prime=j^*F$, in fact they are 
linearly equivalent and  coincide 
outside  $E$. Therefore the scheme-theoretic union $\overline{F}$ of the
scheme-theoretic
images $j_*{\tilde D}^\prime$ and $j_*D^*$ is contained in
$j_*(j^*F)=F$, which is equal to the scheme theoretic union of the strict
images
$j_\#{\tilde D}^\prime=D^\prime$ and $j_\#D^*=D$. Therefore we must have 
$\overline{F}=F$. This implies that $\overline{F}$ cannot have embedded
components, since it is  cut out by a hypersurface on the
arithmetically Cohen-Macaulay scheme $X$.
Since the ideal sheaf of $\overline{F}$ is given by the product of the
ideal sheaves of $j_*(D^*)$ and $D^\prime$, an embedded components of
$j_*(D^*)$ would be an embedded components of $\overline{F}$, but this
is not possible.
On the other side $j_*(D^* + mE)$
for $m>0$ can not be equal to $D$; in fact 
 $D^*+mE\sim (k+1+m)\tilde H -((m+1)f-h-1)\tilde R$
is not contained in $j^*F$ since $j^*F-D^*-mE\sim -mE+(f-h-1)\tilde R$ is 
not effective.
The isomorphism between the ideals follows from
Lemma \ref{ideals}.
\end{pf}

\begin{theorem}
\label{refltwo}
Let $\codim(V, X)=2$. A sheaf $\O_X(a, b)$ with $a, b\in \ZZ$  is
reflexive if and only
if
$b=\deg{\O_{\tilde X}(a, b)}_{|E}<f$.
\end{theorem}
\begin{pf}
The divisors on $\tilde X$ of  type $d\tilde R$ with $d\geq 0$
are proper transforms of divisors on $X$, therefore
the sheaves $\O_X(0, -d)$ are  ideal sheaves on $X$ by Prop. \ref {totale},
and 
therefore
they are reflexive.
The divisors $\sim \tilde H- (f-h-1)\tilde R$ with $0\leq h<f$ are total
transforms of divisors $\sim (h+1)R$ on $X$; by Prop. \ref{totale} 
$\O_X(-1, f-h-1)$ with $0\leq h <f$ are ideal sheaves on $X$, therefore
they are reflexive by Th. \ref{refl}.
By the projection formula 
we have that
$
\O_X(a+a^\prime, b)=
\O_X(a, b)\otimes \O_X(a^\prime, 0)
$
for every $a, a^\prime, b \in \ZZ$.
Since the tensor product of a reflexive sheaf with an invertible sheaf
is reflexive, we obtain that the sheaves $\O_X(a, b)$ are reflexive 
for every $a\in\ZZ$ and $b<f$.
It remains to prove that the sheaves $\O_X(0, d)$ are not reflexive for 
$d\geq f$.
Let us suppose they are reflexive and 
divide $d-1=kf+h$ as usual,
with $k\geq -1$ and $0\leq h <f$;
then the sheaves $\O_X(k, h+1)$ and $\O_X(0, d)$ with
$h<f-1$
(or the sheaves $\O_X(k+1, 0)$ and $\O_X(0, (k+1)f)$  if $h=f-1$)
are both
reflexive and  isomorphic on $X_S$, by Th. \ref{refl}
they should be isomorphic on $X$, but this is not possible.
In fact  the dimension of the respective 
zero  cohomology
groups are different, as one may compute using  (\ref{hi}) and
(\ref{dim}).
\end{pf}

As a consequence of Th. \ref{refltwo} we have that we can write
the divisorial sheaf associated to a Weil divisor $D$ in $X$ in more
than one form. To be more precise, let $d > 0$ and 
 divide $d=kf+h+1\geq 0$, with $k\geq -1$ and $0\leq
h<f$. Let $D\sim dR$  and let us suppose that $D$ is not Cartier, i.e.
  $d\neq 0\mod f$ (otherwise the representation is unique
and it is $\O_X(D)\cong \O_X(k+1, 0)$), then:
\[
\O_X(D)\cong\O_X((k+1),-f+h+1)\cong\O_X(k,h+1).
\]
Let $D\sim -dR$ (now $D$ may be Cartier), then
\[
\O_X(D)\cong\O_X(-(k+1),f-h-1)\cong\O_X(-k, -(h+1))\cong
\cdots\cong\O_X(0, -d).
\]
In the next Corollary we esplicity describe the group $\Div(X)$,
fixing a particular form in which we write a divisorial sheaf.
In
this way it is evident the bijection between
 $\Div(X)$ and
$\Cl(X)\cong \ZZ$.

\begin{corollary}
\label{divtwo}
Let $\codim(V, X)=2$. Then:
\[
\Div(X)=\{\O_X(a, b)\,|\, a, b\in\ZZ, 0\leq b<f \}.
\]
The natural group structure on $\Div (X)$ 
is given by:
$
<\O_X(a, b), \O_X(a^\prime, b^\prime)>\mapsto
\O_X(a+a^\prime+[\frac{b+b^\prime}{f}], b+b^\prime \mod f)
$
and
$
{\O_X(a, b)}^\vee \cong \O_X(-a, -b),
$
where $[\,\,\cdot\,\,]$ denotes the integral part.
Moreover for a fixed $b: 0\leq b <f$, the set
${\Div}_b\{\O_X(a, b)\,|\, a\in\ZZ  \}$ is the set of divisorial sheaves
of Weil divisors $D\sim dR$ with $d = b \mod f$. 
\end{corollary}
\begin{pf}
By Th. \ref{refltwo},
for every $d\in\ZZ$, we can write  the divisorial sheaf associated
to
$D\sim dR$ in the form $\O_X(a, b)$ with $a\in \ZZ$ and $0\leq b<f$. 
\end{pf}

This particular choice is convenient if we want to compute
$h^0(\O_X(D))$, in fact by (\ref{hi}) we know how to
compute it. 
\begin{corollary} 
Let $\codim(V, X)=2$. Let $D\sim dR$ be an effective divisor;
divide $d=kf+h+1$ with $k\geq -1$ and $0\leq h<f$. Let $|D|$ be
the complete linear system of $D$. Then:
\[
\dim |D|= f {{k+r-1}\choose{r}}+(h+2){{k+r-1}\choose{r-1}}-1,
\]
if $h\neq f-1$ (i.e. $D$ is not Cartier), or
\[
\dim |D|= f {{k+r}\choose{r}},
\]
if $D$ is Cartier. 
\end{corollary}
\begin{pf}
Cor. \ref{divtwo}, (\ref{hi}), (\ref{dim}).
\end{pf}

We conclude with the projection formula:
\begin{corollary}{(Projection formula)}
\label{secdiv}
Let $\codim(V, X)=2$. Let $D$ be a Weil divisor on $X$ and let $D^*$ be
its integral total transform.
Then: 
\[
j_*\O_{\tilde X}(D^*)=\O_X(D)
\]
\end{corollary}
\begin{pf}
By Th. \ref{refltwo} $j_*\O_{\tilde X}(D^*)$
 is a reflexive sheaf; since it is the divisorial sheaf associated
to $D$ in the open set $X_S$, then by Th. \ref{refl} 
it is the divisorial sheaf of $D$ an all $X$.
\end{pf}

\sezione{INTERSECTION OF WEIL DIVISORS}

As we have already noted, an effective Weil divisor $D$ on $X\subset\PP^n$
is
a closed
subscheme of $X$ of pure codimension $1$ with no embedded components
(\cite{h2}, Prop. 2.4.). For this reason we can regard $D\subset\PP^n$
as a projective scheme. Given two effective divisors
$D$ and $D^\prime$
on $X$ with no common components, we can consider the 
scheme-theoretic intersection $D\cap D^\prime\subset\PP^n$ and ask for
its degree. 
For {\it degree} of a scheme $Y\subset \PP^n$ we mean the lenght
$h^0(\O_{Y_L})$ of the zero-dimensional scheme $Y_L$ which represents 
a generic ($\codim Y$)-dimensional linear
section of $Y$. As we will see, this problem has an immediate
solution when $\codim(V, X)>2$, via the isomorphism 
$\Cl(X)\cong\CaCl(\tilde X)$. 
In this section we show how to use the integral total transform
to compute this degree in case $\codim(V, X)=2$.

If $\codim(V, X)>2$, by Cor. \ref{Cl} we can define an intersection
form on $\Di(X)$:
\begin{eqnarray*}
I: {\Di(X)}^r &\to &\ZZ \\
<D_1, D_2, \dots, D_r>& \mapsto &
\tilde{D_1}\cdot\tilde{D_2}\cdots\tilde{D_r}
\end{eqnarray*}
exactly as in $\tilde X$; i.e.
$I$  is determined by
the rule:
\[
H^r=f \qquad  H^{r-1}\cdot R=1 \qquad H^{r-2} \cdot R^2=0.
\]
The intersection form $I$ determines 
the degree of the intersection scheme $Y=D_1\cap\cdots \cap D_l\subset
\PP^n$ of $l$ ($l\leq r$)
effective divisors 
which intersect properly.

If $\codim(V, X)=2$, the linear theory of intersection developed
by Mumford in \cite{m} in the case of a normal
 surface can be generalized
to our case by
\[
D_1\cdots D_r=j^*D_1\cdots j^*D_r
\] 
where $j^*D_i$ is the Mumford's total transform of $D_i$.
Given two effective Weil divisors $D$ and $D^\prime$ with no common
components, the 
intersection number
$D\cdot D^\prime\cdot H^{r-2}$ 
does not  represent in this case the degree of the  intersection
scheme $Y=D\cap D^\prime\subset \PP^n$.
To compute this degree we  will find the minimal resolution
of $\O_Y$ as $\O_X$-module (Th. \ref{res}); other applications
of  this resolution  are described 
in the next section. First we need to prove some properties 
of the integral total transform.

\begin{lemma}
\label{genustot}
Let $\codim(V, X)=2$ and let $D$ be an effective
Weil divisor on $X$. Then:
\begin{equation}
\label{genus}
p_a(D)=p_a(D^*).
\end{equation}
\end{lemma}
\begin{pf}
Let $D\sim dR$ and divide $d-1=kf+h$ with $k\geq -1$ and $0\leq h <f$;
since $D^*\sim (k+1)\tilde H-(f-h-1)\tilde R$, with $f-h-1\geq 0$,
from the exact sequence
$
0\to \I_{{D^*}|{\tilde X}}\to \O_{\tilde X}\to \O_{D^*}\to 0
$
by (\ref{Ri}) we  obtain the exact sequence
\[
0\to j_*\I_{{D^*}|{\tilde X}}\to j_*\O_{\tilde X}\to j_*\O_{D^*}\to 0.
\]
Since $j_*\I_{{D^*}|{\tilde X}}=\I_{D|X}$ by Prop. \ref{totale}
and $j_*\O_{\tilde X}=\O_X$,
we get $j_*\O_{D^*}=\O_D$. By (\ref{hi}) we obtain
\[
p_a(D)=1-\chi(\O_D)=1-\chi(\O_{D^*})=p_a(D^*).
\]
\end{pf}

In the next lemma we analize the behaviour of the total transform
with respect to the sum of Weil divisors. First we need the following
definition:

\begin{definition}
With the notation of Def. \ref{deftot}
define: 
\[
\epsilon:=\lceil q\rceil-q.
\]
By definition $\epsilon$
is a rational number in the interval $[0, 1[$ of the kind
$\epsilon=\frac{f-l}{f}$ with
$l=1, 2, \dots, f$. With the notation of (\ref{secdef}) we have: 
\[
\epsilon=\frac{f-h-1}{f}.
\]
We note that given $D\sim dR$ on $X$, $\epsilon$ is uniquely determined by
the class of linear equivalence of $D$, i.e. by $d$.
\end{definition}

\begin{lemma}
\label{sumtot}
Let $\codim(V, X)=2$ and let $D_1\sim d_1R$ and $D_2\sim d_2 R$ be two
effective divisor
on $X$. Then:
\begin{equation}
{(D_1 +D_2)}^*=
\begin{cases}
D_1^*+D_2^*    & \text{if $[\epsilon_1+\epsilon_2]=0$} \\
D_1^*+D_2^* -E & \text{if  $[\epsilon_1+\epsilon_2]=1$} 
\end{cases}
\end{equation}
where $[\,\,\cdot\,\,]$ denotes the integral part.
\end{lemma}
\begin{pf}
By definition \ref{deftot} of integral total transform
we have:
\[
{(D_1+D_2)}^*=\tilde{D_1}+\tilde{D_2}+ \lceil q_1+q_2\rceil E,
\]
since the proper transform of $D_1+D_2$ is $\tilde{D_1}+\tilde{D_2}$.
We find the two cases
$\lceil q_1+q_2\rceil=\lceil q_1\rceil +\lceil q_2\rceil$ if
$0\leq\epsilon_1+\epsilon_2<1$
and $\lceil q_1+q_2\rceil=\lceil q_1\rceil +\lceil q_2\rceil-1$ if
$1\leq\epsilon_1+\epsilon_2<2$.
\end{pf}

\begin{theorem}
\label{reso}
Let $\codim(V, X)=2$ and
let $Y\subset X\subset \PP^n$ be a "complete intersection" of two  
effective
divisors 
$D_1$ and $D_2$  on $X$. The following sequence
is exact:
\begin{equation}
\label{res}
0\to j_*\O_{\tilde X}(-{(D_{1}+D_{2})}^*)
\to 
j_*\O_{\tilde X}(-{D_1}^*)
\oplus j_*\O_{\tilde X}(-{D_2}^*)
\to \I_{Y|X}\to 0
\end{equation}
and therefore it is a "reflexive resolution" of $\O_Y$ as a 
$\O_X$-module.
\end{theorem}
\begin{pf}
First let $f_i:\O_{{D_1}+{D_2}}\to \O_{D_i}$ for $i=1, 2$ be the
projection and let
$g_i:\O_{D_i}\to \O_{{D_1}\cap {D_2}}$. Then there is an exact sequence
of sheaves of $\O_X$-modules:
\begin{equation}
\label{fondseq}
0\to \O_{{D_1}+{D_2}}\stackrel{(f_1,
f_2)}\to\O_{D_1}\oplus\O_{D_2}\stackrel{(g_1,
-g_2)}\to\O_{{D_1}\cap{D_2}}\to 0. 
\end{equation}
For any effective Weil divisor $D$ on $X$ we have already seen the
resolution of $\O_D$:
\begin{equation}
\label{seqone}
0\to j_*\O_{\tilde X}(-D^*)\to\O_X\to \O_{D}\to 0
\end{equation}
in the proof of Lemma \ref{genustot}.
The mapping cone (see \cite{e} pg. 432, pg. 657) between the resolution 
(of type \ref{seqone}) of
$\O_{{D_1}+{D_2}}$
and of $\O_{D_1}\oplus\O_{D_2}$ gives a  resolution
of $\O_{{D_1}\cap{D_2}}$:
\[
j_*\O_{\tilde X}(-{(D_1+D_2)}^*)\hookrightarrow \O_X\oplus j_*\O_{\tilde
X}(-D_1^*)\oplus j_*\O_{\tilde
X}(-D_2^*)\to \O_X\oplus\O_X.
\]
In this resolution we can suppress redundant terms and 
 obtain the required
resolution:
\[
j_*\O_{\tilde X}(-{(D_1+D_2)}^*)\hookrightarrow j_*\O_{\tilde
X}(-D_1^*)\oplus j_*\O_{\tilde
X}(-D_2^*)\to \O_X.
\]
\end{pf}

\begin{note}
The resolution \ref{res} allows us to find a resolution, in general not
minimal, of $\O_Y$ as $\O_{\PP^n}$-module. Indeed, as shown in \cite
{s}, a minimal resolution of
$\O_X(a, b)$ as $\O_{\PP^n}$-module is given by 
 the Eagon-Northcott type complex $C^b(a)$ for $b\geq -1$.
 Since each of the terms in \ref{res}
is a sheaf $\O_X(a, b)$ with $b\geq 0$, then a suitable mapping cone
between the complexes $C^b(a)'s$ gives us the required resolution in
$\PP^n$.
\end{note}

\begin{proposition}
\label{grado}
Let $\codim(V, X)=2$ and  let $D$ and $D^\prime$ be two effective divisors
on $X$ with no common components. Then the degree of the "complete
intersection" scheme $Y=D\cap D^\prime$ is given by:
\begin{equation}
\deg{(D\cap D^\prime)}=
\begin{cases}
D^*\cdot {D^\prime}^*\cdot{\tilde H}^{r-2}    & \text{if
$[\epsilon+\epsilon^\prime]=0$} \\
D^*\cdot {D^\prime}^*\cdot{\tilde H}^{r-2}
 +f(\epsilon+\epsilon^\prime-1)+1 & \text{if
$[\epsilon+\epsilon^\prime]=1$} 
\end{cases}
\end{equation}
\end{proposition}
\begin{pf}
Let us call $X_L$, $Y_L$, $D_L$ and $D^\prime_L$ general
$(r-2)$-dimensional linear sections
of $X$, $Y$,  $D$ and $D^\prime$ respectively, and 
let us call  ${\tilde X}_L$ the canonical
resolution of the rational normal surface $X_L$. 
The resolution of $\O_{Y_L}$
as an $\O_{X_L}$-module is by Th. \ref{reso}:
\[
0\to 
j_*\O_{{\tilde X}_L}(-{(D_{L}+D^\prime_{L})}^*)\to 
j_*\O_{{\tilde X}_L}(-{D_{L}}^*)\oplus 
j_*\O_{{\tilde X}_L}(-{D^\prime_{L}}^*)\to \O_{X_L}. 
\] 
By the proof of Lemma \ref{genustot} we have that
$\chi(j_*\O_{{\tilde X}_L}(-D^*))=p_a(D^*)$ for every effective divisor $D$
on $X_L$.
Let $[\epsilon+\epsilon^\prime]=0$; by Lemma \ref{sumtot}
$({D_{L}+D^\prime_{L})}^*=D_{L}^*+{D^\prime_{L}}^*$. Looking at the
resolution
of $\O_{Y_L}$, by adjunction formula we compute:
\begin{eqnarray*}
h^0(\O_{Y_L}) & =&\frac{1}{2}
(D_{L}^*+{D^\prime_{L}}^*+K_{{\tilde X}_L})\cdot(D_{L}^*+{D^\prime_{L}}^*)\\
& -&\frac{1}{2}(D_{L}^*+K_{{\tilde X}_L})\cdot D_{L}^*
-\frac{1}{2}({D^\prime_{L}}^*+K_{{\tilde X}_L})\cdot {D^\prime_{L}}^*\\
&=& D_{L}^*\cdot {D^\prime_{L}}^* =D^*\cdot {D^\prime}^*\cdot {\tilde
H}^{r-2},
\end{eqnarray*}
where $K_{{\tilde X}_L}\sim -2\tilde H +(f-2)\tilde R$ is the canonical
divisor of ${\tilde X}_L$.
Let $[\epsilon+\epsilon^\prime]=1$; by lemma \ref{sumtot} we have that
$({D_{L}+D^\prime_{L})}^*=D_{L}^*+{D^\prime_{L}}^*-E$.
By an analogous computation we find: 
\begin{eqnarray*}
h^0(\O_{Y_L}) & = &
 D_{L}^* \cdot {D^\prime_{L}}^* -(D_{L}^*+ {D^\prime_{L}}^*)\cdot E
-\frac{1}{2}{(K_{{\tilde X}_L}-E)\cdot E}\\
&= & D_{L}^* \cdot {D^\prime_{L}}^* +(2f-h_1-h_2-2)-f+1\\
&=& D^* \cdot {D^{\prime}}^*\cdot{\tilde H}^{r-2}
+f(\epsilon+\epsilon^{\prime}-1)+1
\end{eqnarray*}
\end{pf}

We note that when $[\epsilon+\epsilon^\prime]=1$, the quantity
$f(\epsilon+\epsilon^\prime-1)$ is bigger than or equal to zero,
therefore in
this case
$\deg(D\cap D^\prime)$ is strictly bigger than the number
$D^*\cdot {D^\prime}^*\cdot{\tilde H}^{r-2}$.
The scheme theoretic intersection $Y$ contains the vertex $V$ as a
component with a certain multiplicity $\geq 0$, we call this number
the {\it integral intersection multiplicity} $m(D, D^\prime; V)$
 of $D$ and $D^\prime$ in $V$, i.e.
\[
m(D, D^\prime; V)= \deg(D\cap D^\prime)-{\tilde D}\cdot {\tilde
D}^\prime 
\cdot{\tilde H}^{r-2}.
\] 
If ${\tilde D}\sim a\tilde H+b\tilde R$
and ${\tilde D}^\prime\sim a^\prime\tilde H+b^\prime\tilde R$
are the proper transforms of $D$ and $D^\prime$, then by Prop.
\ref{grado} we
explicitly compute:
\[
m(D, D^\prime; V)=
\begin{cases}
{\frac{b\cdot b^\prime}{f}}-f(\epsilon\cdot\epsilon^\prime) 
&\text {if $[\epsilon+\epsilon^\prime]=0$}\\
{\frac{b\cdot b^\prime}{f}}-f(\epsilon\cdot\epsilon^\prime)
+f(\epsilon+\epsilon^\prime-1)+1
&\text {if $[\epsilon+\epsilon^\prime]=1$}.
\end{cases}
\]

\begin{note}
The {\it intersection multiplicity}
 of two effective divisors $D$ and $D^\prime$ with no common
components
through the singular locus $V$ on a normal surface $X$, 
is defined in the linear intersection theory of Mumford as the rational
number:
$
i(D, D^\prime; V)
=j^*D\cdot j^*D^\prime - \tilde D\cdot{\tilde
D}^\prime
$.
If $X$ is a rational normal cone, using the same notations as above,
we find that the linear intersection 
multiplicity 
$i(D, D^\prime; V)$ is exactly $\frac{b\cdot b^\prime}{f}$. 
 \end{note}

\sezione{EXAMPLES AND APPLICATIONS}
In this section we show some applications of the previous results. 
In particular in Ex. \ref{quarto}
we use Th. \ref{reso} to compute the arithmetic genus 
of the scheme theoretic intersection $Y$
of two effective divisors on a rational normal cone.

\begin{example}
\label{primo}
In this Example we show that
every effective non degenerate divisor  on a rational normal
cone $X$ of degree
$n-1$
in $\PP^n$ is a curve of maximal arithmetic genus $p_a(C)=G(n, d)$.
\end{example}
Let $C\sim dR$ with $d>n-1$, let us divide $d-1=m(n-1)+\delta$ with $m\geq
1$,
$0\leq\delta\leq n-2$; then $G(n, d)=
{m\choose 2}(n-1)+m\delta$.
By Lemma \ref{genustot} we know that $p_a(C)=p_a(C^*)$, where
$C^*\sim (m+1)\tilde H-(n-2-\delta)\tilde R$. By adjunction formula
on $\tilde X$ we then compute $p_a(C^*)=G(n, d)$. 

\begin{example}
\label{secondo}
Let $\codim(V, X)=2$, then  every effective divisor $D$
and every "complete intersection" $Y$
of two  divisors $D, D^\prime$ on $X\subset \PP^n$ is
arithmetically Cohen-Macaulay.
\end{example}
In the case of one divisor we know by (\ref{seqone}),
(\ref{hi}) and (\ref{cohom})
that
\[
h^i(\I_{D|X}(k))=
h^i(\I_{{D^*}|{\tilde X}}(k))=0
\]
for $1\leq i\leq r-1$ and every $k$. Since $X$ is aritmetically
Cohen-Macaulay from the exact sequence
$
0\to \I_{X|{\PP^n}}\to \I_{D|{\PP^n}}\to \I_{D|X}\to 0
$
we conclude 
$h^i(\I_{D|{\PP^n}}(k))=0$
for $1\leq i\leq r-1$ and every $k$.
Looking at the resolution (\ref{res}) and using (\ref{hi}) and (\ref{cohom})
 we compute:
$
h^i(\I_{Y|X}(k))=0
$
for $1\leq i\leq r-2$ and every $k$; as in the previous case, since  $X$
is aritmetically
Cohen-Macaulay 
we conclude 
$h^i(\I_{Y|{\PP^n}}(k))=0$
for $1\leq i\leq r-2$ and every $k$.

\begin{example}
\label{terzo}
If $\codim(V, X)>2$ and  $Y$ is a "complete intersection" of $l$ ($1\leq 
l\leq r-1$) divisors
$D_i\sim a_i H-b_i R$ with $b_i\geq 0$,
then the resolution of
$\O_Y$
as an $\O_X$-module, 
is a Koszul complex (see \cite{s}, Ex.
3.6.): 
\begin{eqnarray*}
0 & \to &  \O_X(-(a_1+\cdots+a_l), b_1+\cdots+b_l) \to \dots  \\
\dots & \to & 
\sum_{{i_1}<{i_2}} \O_X(-(a_{i_1}+a_{i_2}), b_{i_1}+b_{i_2})
\to  \sum_{i} \O_X(-a_i, b_i) \to \O_X.
\end{eqnarray*}
From this resolution using
(\ref{hi}) and (\ref{cohom})
we conclude that $Y$ is arithmetically Cohen-Macaulay
iff $b_1+\cdots b_l< f$.
\end{example}

\begin{example}
\label{quarto}
Let $\codim(V, X)=2$ and $r\geq 3$.
Resolution \ref{res} can be used to compute  the
arithmetic genus 
of the intersection
scheme $Y$ of two effective divisors $D$ and $D^\prime$ with no common
components. 
\end{example}
Let us suppose that $Y$ is non
degenerate, i.e. $d, d^\prime> f$, then using (\ref{hi}), (\ref{cohom}) 
and (\ref{dim})
we compute:
\begin{eqnarray*}
p_a(Y)&=&
h^0(\O_{\tilde X}(K_{\tilde X}+{(D+ D^\prime)}^*))
-h^0(\O_{\tilde X}(K_{\tilde X}+D^*)-
h^0(\O_{\tilde X}(K_{\tilde X}+{D^\prime}^*)\\
&=& {} \sum_{i=-1}^1 {(-1)}^i ({{\alpha_i-1}\choose{r}}+(f-1-\beta_i)
{{\alpha_i-1}\choose{r-1}})
\end{eqnarray*}
where
$\alpha_{-1}=
\lceil
{\frac {d^{\prime}} {f} }
\rceil$ and $\beta_{-1}=
f
\lceil
{\frac {d^\prime }{f} }
\rceil
-d^\prime
$;
$\alpha_{1}=
\lceil
{\frac {d} {f} }
\rceil$ and $\beta_{1}=
f
\lceil
{\frac {d }{f} }
\rceil
-d
$;
$\alpha_{0}=
\lceil
{\frac {d^{\prime}+d} {f} }
\rceil$ and $\beta_{0}=
f
\lceil
{\frac {d^\prime +d}{f} }
\rceil
-(d^\prime+d)$.
With these notations the degree of $Y$ is:
\[
\deg(Y)=
\sum_{i=-1}^1 {(-1)}^i ({{\alpha_i-1}\choose{r-1}}+(f-1-\beta_i)
{{\alpha_i-1}\choose{r-2}}).
\]

\begin{bibliografia}{XX}


\bibitem{e} D Eisenbud. Commutative Algebra with a View toward
Algebraic Geometry. GTM 150, Springer Verlag, New York, 1994.

\bibitem{ega} A Grothendieck, J Dieudonn\'e. El\'ements de G\'eom\'etrie
Alg\'ebrique. Die Grundlehren der mathematischen Weissenschaften, Band
166. New York: Springer-Verlag, 1971.

 \bibitem{eh} D Eisenbud, J Harris. 
On varieties of minimal degree ( a centennial account), Proceedings of
the
AMS Summer Institute in Algebraic Geometry, Bowdoin, 1985. Proceedings of
Symposia in pure Math., 46, AMS, 1987.


\bibitem{f} R Ferraro: Curve di genere massimo in $ \PP^5$ and
Explicit Resolutions of Double Point Singularities of Surfaces.
PhD dissertation, Universit\`a di Roma "Tor Vergata", 1998.

\bibitem{h1} R Hartshorne. Algebraic Geometry. 
GTM 52, Springer Verlag, New York, 1977.

\bibitem{h2} R  Hartshorne. Generalized divisors on Gorenstein schemes. 
K-Theory Journal, 8, 1994, pp. 287-339.

\bibitem{m} D Mumford. 
The topology of normal singularities of an algebraic surface and a
criterion for simplicity. Publ. Math. I.H.E.S. 9,  1961, pp. 229-246.

\bibitem{r} M Reid. Chapters on Algebraic Surfaces in
Complex Algebraic Geometry. Lectures of a summer programm Park City, 
UT, 1993. Koll\'ar J., IAS/Park City Math. Series 3, AMS, Providence, RI,  
1997, pp. 5-159.

\bibitem{sa} F Sakai.
Weil divisors on normal surfaces.
Duke  Math. J., 51, 4,  1984, pp. 877-887.

\bibitem{s} FO Schreyer. Syzygies of canonical curves and special
linear series. Math. Ann., 275, 1979, pp. 105-137.

\bibitem{se} J. P. Serre. Prolongement de  faiseaux analytique
coherent. Ann. Inst. Fourier, 16, 1966, pp. 363-374.

\end{bibliografia}

\end{document}